\documentclass{article}%
\usepackage{epsfig}
\usepackage{enumerate,graphicx}
\usepackage{amsfonts}
\usepackage{graphicx}
\usepackage{amsmath,amsthm,amssymb}
\usepackage{amssymb}%
\DeclareGraphicsExtensions{.pdf, .jpg, .png , .bmp, .ps}
\usepackage{tikz}
\providecommand{\U}[1]{\protect\rule{.1in}{.1in}}
\newtheorem{thm}{Theorem}[section]
\newtheorem{lem}{Lemma}[section]
\newtheorem{prop}{Proposition}[section]

\theoremstyle{remark}
\newtheorem{rem}[thm]{\bf Remark}

\newcommand\N{{\mathbb {N}}}
\newcommand\R{{\mathbb {R}}}

\newcommand\E{{\mathbb {E}}}

\def\1{{{\mbox{${\rm{1\negthinspace\negthinspace I}}$}}}}

\newcommand{\eref}[1]{(\ref{#1})}

\newcommand\beq{\begin{equation}}
\newcommand\eeq{\end{equation}}

\begin{document}
\title{Deviation inequalities
for separately Lipschitz functionals of composition of random functions}
\author{J\'er\^ome Dedecker,\footnote{Universit\'e Paris Descartes,
Sorbonne Paris Cit\'e,
Laboratoire MAP5
and CNRS UMR 8145, 75016 Paris, France.
E-mail: jerome.dedecker@parisdescartes.fr}\footnotemark[1] \ \
Paul Doukhan\footnote{AGM UMR 8088, Universit\'e Paris Seine, UCP site Saint-Martin,
2 Bd. Adolphe Chauvin 95000 Cergy-Pontoise, France, and CIMFAV, Universidad de Valparaiso, Chile.
Email: doukhan@u-cergy.fr}\footnotemark[2]   \ \   and   \ \  Xiequan Fan\footnote{Center for Applied Mathematics,
Tianjin University, Tianjin, China.
E-mail: fanxiequan@hotmail.com}\footnotemark[3] \footnote{All authors contributed equally to this work}}

\date{}

\maketitle

\abstract{We consider a class of non-homogeneous Markov chains, that contains many natural examples. Next,  using martingale methods, we establish some  deviation and moment  inequalities for separately Lipschitz functions of such a chain,  under  moment conditions on some dominating random variables.}

\medskip


\medskip


\section{Introduction}
Concentration inequalities are essential tools for ensuring the validity of many statistical procedures; let us cite for instance \cite{FMS} for classification problems,
\cite{M03} for model selection, and \cite{BRT} for high dimensional procedures (see also \cite{AD11} and \cite{Z17} in a dependent framework).

These inequalities are not easy to prove in a dependent context; up to now this has  been done under quite restrictive assumptions,
and mainly for bounded functionals of the variables in a stationary/homogeneous context. As a non exhaustive list,
let us quote \cite{Sa00}, \cite{Ri00}, \cite{DGW04}, \cite{Ad08}, \cite{DF15} and \cite{DG15}.
Among these references, the case of unbounded functionals has been  investigated in \cite{Ad08}  for geometrically ergodic Markov chains,  and in
\cite{DF15} for iterated random functions satisfying a mean-contraction condition (see condition \eqref{contract} below, with $F_n=F$).

In the paper \cite{DF15}, the authors obtained very precise inequalities for Lipschitz
functionals of the chain, by assuming moreover a Lipschitz condition on the function that generates the chain (see condition  \eqref{c2old} below).
However, this last condition is in fact quite restrictive, and does not hold for many natural models satisfying the  mean-contraction  property.

In the present paper we  enlarge the class of Markov chains studied in \cite{DF15}, by considering non-homogeneous Markov chains obtained through composition of random functions
(see the recursive mechanism \eqref{Mchain} below), and by making no extra assumptions  than the (uniform) mean-contraction \eqref{contract}.
As in \cite{DF15},
we shall use the decomposition of the functional of the chain in terms of martingale differences, as first introduced by Yurinskii \cite{Yu}. This method  is well adapted to the non-homogeneous Markov context, because it is   intrinsically a non-stationary method.
Following this approach, we obtain deviation and moment inequalities for separately Lipschitz functionals of the chain that are driven by the distribution of some
dominating random variables.

The present article was conceived within the general framework  of non stationary time series, which is now widely discussed in the context global warming \cite{C00}.
Besides temperatures or ozone concentration, most of the real life phenomena present trends and periodicities. A first excellent view of those questions may be found in \cite{BD87}, but this is a  linear view of time series analysis. It appears important to consider cases where the dynamic of the models itself is non time-homogeneous.
For instance  \cite{Da18} and \cite{BD17} provide different views for a more relevant  dynamical approach including local stationarity and non-periodic features.
The present paper aims at developing reasonable concentration and probability inequalities for  non-homogeneous Markov chains able to model some of the above features.

Before going into further details, let us give a simple class of examples to which our results apply.
We consider a generalized ${\mathbb R}^k$-valued  auto-regressive processes
\begin{equation}\label{GARM}
X_{n}= A_n X_{n-1} + B_n  \, ,
\end{equation}
where $A_n$ is a random $k \times k$  matrix and $B_n$ is an ${\mathbb R}^k$-valued random variable. Here $\varepsilon_n= (A_n, B_n)$ are independent random
variables, and $(\varepsilon_n)_{n \geq 2}$ is independent of the initial random variable $X_1$. Let $| \cdot |$ be a norm on ${\mathbb R}^k$.
Then, the Markov chain $X_n$  satisfies the mean contraction condition  \eqref{contract} for the norm  $| \cdot |$  as soon as
\begin{equation}\label{MCGA}
\sup_{n \geq 2} {\mathbb E}[|A_n|] \leq \rho \quad \text{for some $\rho <1$},
\end{equation}
where as usual $|A_n|= \sup_{|x|=1} |A_n x|$.

Model \eqref{GARM} contains a lot of natural examples (see for instance \cite{DF}, Sections 2.2 and 2.3), but does not fit within the framework of \cite{DF15}; moreover, it
has no reason to be mixing in the sense of Rosenblatt \cite{Ro} without further assumptions on the distribution of $(\varepsilon_n)_{n \geq 2}$. Recall that the chain
$X_n$ is non-homegeneous, since we do not assume here that the variables $\varepsilon_n$ are iid.

Let now $f: ({{\mathbb R}^k})^n \mapsto {\mathbb R}$ be a separately Lipschitz function, such that
\begin{equation} \label{SepLip}
|f(x_1, x_2, \ldots, x_n)-f(x'_1, x'_2, \ldots, x'_n)| \, \leq \, |x_1-x'_1|+\cdots +  |x_n-x'_n| \, ,
\end{equation}
and let also $S_n= f(X_1, \ldots, X_n) - {\mathbb E} [f(X_1, \ldots, X_n)]$.

For simplicity, let us consider the case where the chain starts at $X_1=0$. Assuming  that
$ \|A_n\|_p^p:= {\mathbb E}[|A_n|^p] < \infty$ and $\| B_n\|_p^p:={\mathbb E}[|B_n|^p]< \infty$ for any $n\geq 2$ and some $p>1$, we infer from \eqref{HGA}
(control of the ${\mathbb L}^p$-norm of the dominating  variables $H_k(X_{k-1}, \varepsilon_k)$ defined in \eqref{defHk}) and Propositions \ref{MZIP} and \ref{VBEI} that
$$
\| S_n \|_p^p \leq    C(p, \rho)  \left(
 \sum_{k=3}^n \|X_{k-1}\|_p^2 \|A_k\|_p^2 + \sum_{k=2}^n   \|B_k\|_p^2 \right)^{2}  \quad \text{if $p\geq 2$},
$$
and
$$
\| S_n \|_p^p \leq    C(p, \rho)  \left (
 \sum_{k=3}^n \|X_{k-1}\|_p^p \|A_k\|_p^p + \sum_{k=2}^n   \|B_k\|_p^p \right ) \quad \text{if $p \in (1,2)$},
$$
for some constant $C(p, \rho)$ depending only on $(p, \rho)$. These inequalities are satisfactory, because if $X_n=B_n$ for $n \geq 2$ (case $A_n =0$), we recover  for $p \geq 2$ the
usual Marcinkiewicz-Zygmund inequalities  (see \cite{R09}) for ${\mathbb L}_p$-norms of sums of independent random variables, and for $p \in (1,2)$ the usual von Bahr-Esseen inequalities
(see \cite{P10}).  Note that, under the stronger condition  than \eqref{MCGA}: $\sup_{n \geq 2} \|A_n\|_p \leq \rho$
(${\mathbb L}^p$-contraction), and if $\sup_{n \geq 2} \|B_n\|_p  < \infty$, we obtain that $\|S_n\|_p = O(\sqrt n)$ if $p\geq 2$ and $\|S_n\|_p = n^{1/p}$ if $p \in (1,2)$, which is exactly what we could expect
for ${\mathbb L}^p$-norms of partial sums in  a quasi-stationary regime.

Under more restrictive conditions on $(A_n, B_n)_{n \geq 2}$, one can also obtain some semi-exponential bounds for the deviation
of $S_n$. For the sake of simplicity, let us assume  that  $\sup_{n \geq 2} \|A_n\|_\infty \leq \rho$ (uniform contraction) and that
there exist $\kappa >0, \alpha \in (0,1)$ such that
\begin{equation}\label{unifexp}
  \sup_{k \geq 2}  {\mathbb E}\left [ \exp   \left \{  \kappa |B_k|^{\frac{2\alpha}{1- \alpha}} \right \}\right ] < \infty  \, .
\end{equation}
It is then easy to see that dominating variables $H_k(X_{k-1}, \varepsilon_k)$ defined in \eqref{defHk} also satisfy the uniform bound
\eqref{unifexp} (for the same $\alpha$ and a different $\kappa$, say $\kappa'$). Hence, it follows from Proposition \ref{propExp} that
\begin{eqnarray}\label{finintro}
  {\mathbb P}\left( | S_n | \geq  n x\right)
   \leq   C( x)  \exp\left\{- K x^{2\alpha} n^\alpha \right\} \, ,
\end{eqnarray}
for any $x>0$,
where  the positive constant $K$ depends only on $(\rho, \alpha, \kappa')$, and
$$
C(x)=  2+ c(\alpha, \rho, \kappa')    \left( \frac{1 }
{ x^{2\alpha} }  + \frac{1}{  x^2 } \right) \, .
$$
In particular, we obtain from \eqref{finintro} the following moderate deviation behavior: for any $\delta \in (1/2, 1]$, there exist $a>0$ such that
\begin{equation}\label{etlabete}
 {\mathbb P}\left( | S_n | \geq  n^\delta \right)= O \left ( \exp \left \{ - a n^{\alpha(2\delta -1)}  \right \} \right ) \, .
\end{equation}
 Note that, for $\delta=1$,  this is in accordance with the best possible
rate for large deviation of sums of martingale differences (see Theorem 2.1 in \cite{Fx1}).

\section{Composition of random functions}\label{SEC2}

Let $(\Omega, {\mathcal A}, {\mathbb P})$ be a probability space.
Let $({\mathcal X}, d)$ and  $({\mathcal Y}, \delta)$ be two complete separable metric spaces.
Let $(\varepsilon_i)_{i \geq 2}$ be a sequence of independent  ${\mathcal Y}$-valued
random variables. Let $X_1$ be a ${\mathcal X}$-valued random variable independent of $(\varepsilon_i)_{i \geq 2}$. We consider the Markov chain $(X_i)_{i \geq 1}$
such that
\begin{equation}\label{Mchain}
X_n=F_n(X_{n-1}, \varepsilon_n), \quad \text{  $n\geq 2$},
\end{equation}
where $F_n: {\mathcal X} \times {\mathcal Y} \rightarrow {\mathcal X}$
is such that
\begin{equation}\label{contract}
{\mathbb E}\big [ d\big(F_n(x, \varepsilon_n), F_n(x', \varepsilon_n)\big) \big] \leq
\rho \, d(x, x')
\end{equation}
for some constant $\rho \in [0,1)$ not depending on $n$.

\smallskip

In the paper \cite{DF15}, the authors studied a class of homogeneous Markov chains (that is, with  $F_n=F$ and   $(\varepsilon_i)_{i \geq 2}$ a sequence of  i.i.d.\ random variables) satisfying \eqref{contract} and
the condition
\begin{equation}\label{c2old}
d(F(x,y), F(x,y')) \leq C \,\delta(y,y')
\end{equation}
for some positive constant $C$. Under this additional constraint, they obtained  very precise upper bounds for the deviation of separately Lipschitz functionals
of the chain; this is possible, because in that case, the martingales differences $M_k$ from McDiarmid's decomposition are bounded by a function of $\varepsilon_k$, which is then independent of the past $\sigma$-field of the chain.

However, condition \eqref{c2old} is quite restrictive, and is not satisfied for many natural models (a short list of such models is presented below). In the present paper, we shall not assume that
\eqref{c2old} is satisfied. In this more general setting, the dominating random variables are
\begin{equation}\label{defHk}
H_k(X_{k-1}, \varepsilon_k) \quad \text{where} \quad
H_{k}(x,y)= \int d(F_k(x,y),F_k(x,y'))P_{\varepsilon_k}(dy')
\end{equation}
(see Proposition \ref{McD} below).  The main difference with \cite{DF15} is that these dominating random variables are no longer independent from the past
$\sigma$-field of the chain. Hence, the  deviations bounds that we obtain are not as precise  as in \cite{DF15}, but apply to a much larger class of
(non homogeneous) Markov chains.

\smallskip

\begin{rem}\label{remalpha}
Note that if \eqref{contract} holds for the distance $d$, then, for any $\alpha \in (0,1]$,  it also hold for the distance
$d_\alpha(x,y)=(d(x,y))^\alpha$ with $\rho^\alpha$ instead of $\rho$.
 This is elementary, but nevertheless important:
it means that we can also obtain concentration inequalities for separately Lipschitz functions with respect to $d_\alpha$ by controlling
the behavior of $H_{k,\alpha}(X_{k-1}, \varepsilon_k)$ (whose definition is as in \eqref{defHk} for the distance $d_\alpha$).
Note that separately Lipschitz functions with respect to $d_\alpha$ are less and less regular as $\alpha$ approaches 0.

\end{rem}

\smallskip

\begin{rem}
Let us quote an error in the paper \cite{DF15}. The inequality (1.4) of that paper gives an upper bound for the
quantity ${\mathbb E}[H(d(X_n,x_0))]$ when $H$ is any increasing function from ${\mathbb R}^+$ to ${\mathbb R}^+$.
However this upper bound is not true in general under the assumption (1.2) of \cite{DF15} (which is similar to our assumption
\eqref{contract}), but it holds  under the much more restrictive assumption $  d(F(x, y), F(x', y))  \leq
\rho \, d(x, x')$. The error comes from the fact that the first version of the paper  \cite{DF15} was written under this more restrictive assumption.
Note that this wrong inequality was not used at any points in the proofs of the main results in \cite{DF15}, but only in Items 4 of Remarks 3.1 and 3.2
(which are therefore not correct).

\end{rem}

\subsection{Examples}
In this subsection, we give a non exhaustive list of models satisfying condition \eqref{contract}, and we show how to control the
moments of the dominating  variables $H_k(X_{k-1}, \varepsilon_k)$ defined by \eqref{defHk}. For the sake of simplicity, we shall only deal with the moments
of order $p$ of $H_k(X_{k-1}, \varepsilon_k)$, but similar computations may be done for exponential moments.
We refer to \cite{DFT12} and \cite{Da18} for more examples.

\begin{itemize}
\item
{\bf ARCH-type  models}. For which
$$
F_n(x,y)=M_{\theta_n}(x,y), \qquad \mbox{with}\quad M_{\theta}(x,y)=\sqrt{a^2x^2+b^2} \cdot y, \quad \theta=(a,b)\, .
$$
Thus the non stationarity appears simply from changes in the parameter $\theta_n=(a_n,b_n)$.
In that case \eqref{contract} is satisfied for $d(x,x')= |x-x'|$ provided $\sup_{n \geq 2}|a_n|{\mathbb E}[ |\varepsilon_n| ] = \rho$
for some $\rho <1$.

\medskip

For these models
$$
H_k(X_{k-1}, \varepsilon_k)= \sqrt{a_k^2X_{k-1}^2+b_k^2} \int | \varepsilon_k-y| P_{\varepsilon_k}(dy) \, ,
$$
and the moments of order $p$ of $H_k(X_{k-1}, \varepsilon_k)$ satisfy
$$
  \| H_k(X_{k-1}, \varepsilon_k) \|_p^p \leq 2^{p-1} {\mathbb E} \left [ \left (a_k^2X_{k-1}^2+b_k^2 \right )^{p/2}\right ] \| \varepsilon_k \|_p^p \, , \ \ \ p\geq 1.
$$

Those models are easy to extend in an $\R^k$-valued framework.
For instance, one can consider
$X_{n}= A_n(X_{n-1}) \, \varepsilon_n$, where $A_n(x)$ is a $k\times k$ matrix and $\varepsilon_n$  are ${\mathbb R}^k$-valued random variables.
Let $| \cdot |$ be a norm on ${\mathbb R}^k$, and $|A|= \sup_{|x|=1} |A_n x|$ be the associated
  matrix norm.  Now, if
$|A_n(x) -A_n(x')|  \leq a_n |x-x'|$,
then the  condition \eqref{contract} is satisfied as soon as  $\sup_{n \geq 2} a_n {\mathbb E}[ |\varepsilon_n| ] = \rho$
for some $\rho <1$.

\item  {\bf Switching models.} Many analogous models can be provided with a switching, e.g.\ for the first ARCH-model, such  a parametric model is given with ${\cal X}=\R$, ${\cal Y}=\R\times\{0,1\}$, a parameter $\theta=(a,b,a',b')\in \R^4$ and
$$
 M_\theta(x,y)=y_2\sqrt{a^2x^2+b^2} \cdot y_1+(1-y_2)\sqrt{a'^2x^2+b'^2} \cdot y_1 \, .
$$
Here $(\varepsilon_n)_{n \geq 2}$ is a sequence of independent random variables with values in ${\mathbb R} \times \{0, 1 \}$. Using the notation
$\varepsilon_n=(\varepsilon_n^{(1)}, \varepsilon_n^{(2)})$, we see that condition \eqref{contract} is satisfied as soon as
$$
  \sup_{n \geq 2} \  \left ( (1-\E [\varepsilon_n^{(2)}])|a_n|+ \E[\varepsilon_n^{(2)}]|a_n'| \right)\E[|\varepsilon_n^{(1)}|]= \rho
$$
for some $\rho <1$.

Now, similar computations as for the first example lead to
\begin{multline*}
  \| H_k(X_{k-1}, \varepsilon_k) \|_p^p \leq 4^{p-1} {\mathbb E} \left [ \left (a_k^2X_{k-1}^2+b_k^2 \right )^{p/2}\right ] \left \| \varepsilon_k^{(1)}  \varepsilon_k^{(2)} \right \|_p^p \\
  + 4^{p-1} {\mathbb E} \left [\left (a'^2_kX_{k-1}^2+b'^2_k \right )^{p/2}\right ]\left \| \varepsilon_k^{(1)} (1- \varepsilon_k^{(2)}) \right \|_p^p , \, \ \ \ p\geq 1.
\end{multline*}

\item {\bf Generalized ${\mathbb R}^k$-valued  auto-regressive processes}. We consider here the Model \eqref{GARM} presented in the introduction.
Recall that $A_n$ is a random $k \times k$  matrix and $B_n$ is an ${\mathbb R}^k$-valued random variable. Here $\varepsilon_n= (A_n, B_n)$ are independent random
variables, and $(\varepsilon_n)_{n \geq 2}$ is independent of the initial random variable $X_1$.  Model \eqref{GARM} is a composition of random functions as in \eqref{Mchain},
with
$$F_n(x,y)=F(x,y)=y_1 x + y_2 \, . $$

Let $| \cdot |$ be a norm on ${\mathbb R}^k$, and let  as usual $|A_n|= \sup_{|x|=1} |A_n x|$. The condition \eqref{contract} is satisfied as soon as
\eqref{MCGA} holds.

\medskip

For these models
$$
H_k(X_{k-1}, \varepsilon_k) \leq | X_{k-1}|  \int  |A_k-y  | P_{A_k}(dy) +  \int | B_k-y| P_{B_k}(dy) \, ,
$$
and the moments of order $p$ of $H_k(X_{k-1}, \varepsilon_k)$ satisfy
\begin{equation}\label{HGA}
  \| H_k(X_{k-1}, \varepsilon_k) \|_p^p \leq 4^{p-1} {\mathbb E} \left  [ |X_{k-1}|^p \right ] {\mathbb E} \left [ |A_k|^p \right ] + 4^{p-1} {\mathbb E} \left [ |B_k|^p \right ] , \ \ \ \ p\geq 1.
\end{equation}

\item {\bf INAR($1$) type models.} In this case, let $y=(y_0,y_1,y_2,\dots,y_p,\ldots)\in {\cal Y}=\N^{\N}$ and   $\varepsilon=(\varepsilon^{(0)},\varepsilon^{(1)},\ldots,\varepsilon^{(p)},\ldots)$, where $(\varepsilon^{(1)},\ldots,\varepsilon^{(p)},\ldots)$ is a sequence of i.i.d. integer valued random variables. The function $F$ is then given by
 $$
F(x,y)=y_0+\textbf{1}_{\{x\neq 0\}}\sum_{k=1}^x y_k\,.
$$
Here,  $(\varepsilon_n)_{n \geq 2}$ is an i.i.d.\ sequence distributed as $\varepsilon$. It is then easy to see that \eqref{contract} is satisfied
provided $\rho={\mathbb E}[\varepsilon^{(1)}]< 1$.

\medskip

We shall now give some hints to control the moments of the dominating random variables  $H_k(X_{k-1}, \varepsilon_k)$.
Let $ \tilde \varepsilon_k$ be distributed as $\varepsilon_k$ and  independent of $(\varepsilon_k, X_{k-1})$. We then have
that
$$
 H_k(X_{k-1}, \varepsilon_k) = {\mathbb E} \big[d(F_k(X_{k-1},\varepsilon_k),F_k(X_{k-1},\tilde \varepsilon_k)) \big| X_{k-1}, \varepsilon_k \big]   .
$$

For the INAR(1) model, we have
$$
H_k(X_{k-1}, \varepsilon_k)= \E \bigg[ \Big | (\varepsilon_k^{(0)}- \tilde \varepsilon_k^{(0)})+\textbf{1}_{\{X_{k-1}\neq 0\}}\sum_{i=1}^{X_{k-1}}
 ( \varepsilon^{(i)}_k - \tilde \varepsilon^{(i)}_k)\Big |  \  \bigg  | X_{k-1}, \varepsilon_k \bigg]\\
$$
By contraction, we get that
$$
 \| H_k(X_{k-1}, \varepsilon_k) \|_p^p \leq  \bigg \| (\varepsilon_k^{(0)}- \tilde \varepsilon_k^{(0)})+\textbf{1}_{\{X_{k-1}\neq 0\}}\sum_{i=1}^{X_{k-1}}
 ( \varepsilon^{(i)}_k - \tilde \varepsilon^{(i)}_k)\bigg \|_p^p , \, \ \ \ \ p\geq 1.
$$
For $p \geq 2$, applying the Marcinkiewicz-Zygmund inequality given in \cite{R09}, we  get that
$$
 \| H_k(X_{k-1}, \varepsilon_k) \|_p^p \leq  2^{(p-2)/2}  \|  \varepsilon^{(0)}-\tilde \varepsilon^{(0)}\|_p^p +  (p-1)^{p/2} 2^{(p-2)/2} {\mathbb E} \left [X_{k-1}^{ p/2} \right] \|  \varepsilon^{(1)}-\tilde \varepsilon^{(1)}\|_p^p \, .
$$
For $p \in (1,2)$, applying  the von-Bahr Essen inequality given in \cite{P10}, we get that
$$
 \| H_k(X_{k-1}, \varepsilon_k) \|_p^p \leq    \|  \varepsilon^{(0)}-\tilde \varepsilon^{(0)}\|_p^p +  2^{2- p} {\mathbb E} \left[X_{k-1} \right] \|  \varepsilon^{(1)}-\tilde \varepsilon^{(1)}\|_p^p \, .
$$

Note that non-stationary variants of this model can be obtained by considering  independent (but non i.i.d.)\ $\varepsilon_n$'s, with the constraint:
$\sup_{n\geq 2} {\mathbb E}[\varepsilon_n^{(1)}] < 1$.
\item
{\bf GLM-Poisson models.} Besides the standard ARCH-models the simplest case is that of Poisson ARCH--models, where $(\varepsilon_n)_{n \geq 2}$ is a
sequence of i.i.d.\ unit Poisson processes.
Consider a sequence of  functions $f_n:\mathbb{N}\to \mathbb{R}^+$
and set
$$
F_n(x,y)=y(f_n(x)),
$$
where $y:\mathbb{R}^+\to \mathbb{N}$ denotes a function.
In that case, the condition \eqref{contract} is satisfied if $|f_n(x)-f_n(x')|\leq \rho |x-x'|$ for any $n\geq 2$, any $x,x' \in {\mathbb N}$, and some $\rho <1$.

\medskip

For these models,
$$
H_k(X_{k-1},\varepsilon_k)= \E\big[| \varepsilon(f_k(X_{k-1}))-\varepsilon_k(f_k(X_{k-1}))|\big|\varepsilon_k,X_{k-1}\big]   ,
$$
where $ \varepsilon$ is a unit Poisson process independent of $(X_{k-1}, \varepsilon_k)$. By contraction, we get that
$$
 \| H_k(X_{k-1}, \varepsilon_k) \|_p^p \leq 2^{p-1} \| \varepsilon_k(f_k(X_{k-1}))\|_p^p= 2^{p-1}  {\mathbb E} \left [Q_p(f_k(X_{k-1})) \right ]   ,
$$
where $Q_p(t)=\| \varepsilon(t)\|_p^p$. Note that, when $p$ is an integer, $Q_p$  denotes the Stirling polynomial defined through Stirling numbers
(see the Lemma A-1 in \cite{DFT12} from \cite{FLO06}).
\item {\bf GLM--GARCH Poisson models}.
One can give numerous extensions of the previous model. Keeping the same notations, one can consider
$f_n:  \mathbb{R}^+ \times \mathbb{N}\to \mathbb{R}^+$,
and
$$
F_n(x,y)=(f_n(x), y(f_n(x)),
$$
where $x=(\lambda, z) \in  \mathbb{R}^+ \times \mathbb{N}$.
Let
$|x|=|\lambda|+a|z|$. Then
$$
{\mathbb E} [ | F_n(x, \varepsilon)- F_n(x', \varepsilon_n)|] \leq (1+a) | |f_n(x)-f_n(x')| \, ,
$$
and
\eqref{contract} is true provided that
$|f_n(x)-f_n(x')|\leq L_n |x-x'|$ and $\sup_{n \geq 2} L_n (1+a) = \rho $ for some $\rho <1$.

\medskip

For these models,
$$
H_k(X_{k-1},\varepsilon_k)= a \E\big[| \varepsilon(f_k(X_{k-1}))-\varepsilon_k(f_k(X_{k-1}))|\big|\varepsilon_k,X_{k-1}\big]   ,
$$
where $ \varepsilon$ is a unit Poisson process independent of $(X_{k-1}, \varepsilon_k)$. Hence, the moments of order $p$
of $H_k(X_{k-1},\varepsilon_k)$ can be controlled exactly as in the previous example.
\end{itemize}

\section{Separately Lipschitz functions of $X_1, \ldots , X_n$}
\setcounter{equation}{0}
Let $f: {\mathcal X}^n \mapsto {\mathbb R}$ be a separately Lipschitz function, such that
\begin{equation} \label{codiMD}
|f(x_1, x_2, \ldots, x_n)-f(x'_1, x'_2, \ldots, x'_n)| \, \leq \, d(x_1,x'_1)+\cdots +  d(x_n, x'_n) \, .
\end{equation}
Let
\begin{equation}\label{Sn}
S_n \, := \, f(X_1, \ldots , X_n) -{\mathbb E}[f(X_1, \ldots , X_n)]\, .
\end{equation}
We   introduce the natural filtration of the chain, that is  ${\mathcal F}_0=\{\emptyset, \Omega \}$
and for all $k \in {\mathbb{N}}^{*}$,
${\mathcal F}_k= \sigma(X_1, X_2,  \ldots, X_k)$.
Define
\begin{equation}\label{gk}
g_k(X_1, \ldots , X_k) \, = \, {\mathbb E}[f(X_1, \ldots, X_n)|{\mathcal F}_k]
\end{equation}
and
\begin{equation}\label{dk}
M_k\,= \, g_k(X_1, \ldots, X_k)-g_{k-1}(X_1, \ldots, X_{k-1})  .
\end{equation}
For all $k \in [1, n-1]$, let
$$
S_k:=M_1+M_2+\cdots + M_k  ,
$$
and notice that, by the definition of $M_k$'s, the functional $S_n$ introduced in
\eref{Sn} satisfies
$$
S_n \, = \, M_1+M_2+\cdots + M_n   .
$$
Thus $S_k$ is a martingale adapted to the natural filtration ${\mathcal F}_k$.
This representation appears in
Yurinskii \cite{Yu} and in p.\ 33 of the monograph by Milman and Schechtman \cite{MS}.
In the setting of separately Lipschitz functions of independent random variables ({\it i.e.}\ when $X_i=\varepsilon_i$)
it has been used by  McDiarmid \cite{M}  to get an exponential bound
on tail probabilities ${\mathbb P}(S_n \geq x), x\geq 0$.

The following Proposition, similar to Proposition 2.1 in \cite{DF15},  collects some interesting properties of
the function  $g_k$ and of the martingale difference $M_k$.

\begin{prop}\label{McD}
For all $k \in {\mathbb N}$
and any $\rho$ in $[0, 1)$, denote $K_k(\rho)=(1-\rho^{k+1})/(1-\rho)= 1 + \rho + \cdots + \rho^{k}$.
Let $(X_i)_{i \geq 1}$ be a Markov chain satisfying \eref{Mchain} for some
functions $F_n$ satisfying \eref{contract}. Let $g_k$ and $M_k$ be defined
by
\eref{gk} and \eref{dk} respectively.
\begin{enumerate}
\item
The function $g_k$ is separately Lipschitz, and satisfies
$$
\big|g_k(x_1, x_2, \ldots, x_k)-g_k(x'_1, x'_2, \ldots, x'_k) \big| \leq d(x_1,x'_1)+\cdots + d(x_{k-1}, x'_{k-1})+ K_{n-k}(\rho) d(x_k, x'_k) \, .
$$
\item Denote by $P_{X_1}$  and $P_{\varepsilon_k}$   the distribution of $X_1$ and the
distribution of  $\varepsilon_k$  respectively. Let $G_{X_1}$ and $H_{k}$ be  two functions defined as follows
$$
G_{X_1}(x)=\int d(x, x') P_{X_1}(dx') \quad \text{and}
\quad
H_{k}(x,y)= \int d(F_k(x,y),F_k(x,y'))P_{\varepsilon_k}(dy') \, .
$$
Then, the martingale difference  $M_k$  satisfies
$$
   |M_1| \leq K_{n-1}(\rho) G_{X_1}(X_1) \quad
   \text{and } \quad
   |M_k|\leq K_{n-k}(\rho)  H_{k}(X_{k-1}, \varepsilon_k), \ \  \, k \in [2, n].
$$
\end{enumerate}
\end{prop}
\begin{rem}\label{remc2}
Assume moreover that $F_n$ satisfies
\begin{equation}\label{c2}
d(F_n(x,y), F_n(x,y')) \leq C(x) \,\delta(y,y')
\end{equation}
for some function $C(x)  \geq 0$ not depending on $n$, and
let $G_{k}$ be the  function defined by
$$
G_{k}(y)= \int\delta(y,y')P_{\varepsilon_k}(dy') \, \
   \text{ for all $k \in [2, n]$}.
$$
Then $H_{k}(x,y) \leq C(x)G_{k}(y)$
and, consequently,
$$
 |M_k|\leq K_{n-k}(\rho)C(X_{k-1})  G_{k}(\varepsilon_k) \, \
   \text{ for all $k \in [2, n]$}.
$$
Note that \eqref{c2} is a non-uniform version of \eqref{c2old}, which is satisfied for many examples (for instance the three first
examples of Section 2). However, it seems quite difficult to check for INAR or GLM type models, while the moments of the dominating variables
$H_k(X_{k-1}, \varepsilon_k)$ are easy to control for such models (see Section \ref{SEC2}).
\end{rem}

\noindent {\emph {Proof}.} The first point will be proved by recurrence in the backward sense.
For $k=n,$ the result is obvious due to  $g_n=f$. Suppose it is true at step $k$, and let us prove it
at step $k-1$. By definition
$$
g_{k-1}(X_1, \ldots, X_{k-1})={\mathbb E}[g_k(X_1, \ldots, X_k)|{\mathcal F}_{k-1}]= \int g_k(X_k, \ldots, X_{k-1}, F_k(X_{k-1},y)) P_{\varepsilon_k}(dy)\, .
$$
Then it is easy to see that
\begin{multline}\label{triv1}
|g_{k-1}(x_1, x_2, \ldots, x_{k-1})-g_{k-1}(x'_1, x'_2, \ldots, x'_{k-1})|  \\ \leq \, \int \Big|g_k(x_1, x_2, \ldots, F_k(x_{k-1},y))-g_k(x'_1, x'_2, \ldots, F_k(x'_{k-1},y)) \Big|
P_{\varepsilon_k}(dy)\, .
\end{multline}
Now, by assumption and condition \eref{contract},
\begin{multline}\label{triv2}
\int \Big|g_k(x_1, x_2, \ldots, F_k(x_{k-1},y))-g_k(x'_1, x'_2, \ldots, F_k(x'_{k-1},y)) \Big|
P_{\varepsilon_k}(dy)\\
\leq d(x_1,x'_1)+\cdots + d(x_{k-1}, x'_{k-1})+ K_{n-k}(\rho) \int d(F_k(x_{k-1},y), F_k(x'_{k-1},y))P_{\varepsilon_k}(dy)
\\
\leq d(x_1,x'_1)+\cdots  + (1+\rho K_{n-k}(\rho)) d(x_{k-1}, x'_{k-1}) \\
\leq d(x_1,x'_1)+\cdots  +  K_{n-k+1}(\rho) d(x_{k-1}, x'_{k-1}) \, .
\end{multline}
The point 1 follows from \eref{triv1} and \eref{triv2}.

Next, we  prove the point 2. First notice that
$$
|M_1|=\bigg|g_1(X_1)-\int g_1(x)P_{X_1}(dx)\bigg|\leq
 K_{n-1}(\rho)\int d(X_1,x)P_{X_1}(dx)=K_{n-1}(\rho) G_{X_1}(X_1) \, .
$$
Similarly, for all $k \geq 2$,
\begin{align*}
 |M_k| \,& =\, \Big|g_k(X_1, \cdots, X_k)-{\mathbb E}[g_k(X_1, \cdots, X_k)|{\mathcal F}_{k-1}]\Big|\\
 &\leq \, \int \Big |g_k(X_1, \cdots, F_k(X_{k-1}, \varepsilon_k))-g_k(X_1, \cdots, F_k(X_{k-1}, y))\Big| P_{\varepsilon_k}(dy)\\
 &\leq \, K_{n-k}(\rho)\int  d(F_k(X_{k-1},\varepsilon_k),
 F_k(X_{k-1}, y))P_{\varepsilon_k}(dy)
 \, = \, K_{n-k}(\rho)
  H_{k}(X_{k-1},\varepsilon_k) \, .
\end{align*}
This completes the proof of Proposition \ref{McD}.
\hfill\qed

\section{Deviation inequalities for the functional $S_n$} \label{deviationiq}

\setcounter{equation}{0}
Let $(X_i)_{i \geq 1}$ be a Markov chain satisfying \eref{Mchain} for some
functions $F_n$ satisfying \eref{contract}.
In this section, we apply inequalities for martingales to bound up
the deviation of the functional $S_n$ defined by \eref{Sn}. Some of these inequalities
are direct applications of known inequalities, and some deserve a short proof.

Denote by $S_{2,n}= S_{n}-M_1$, and let $a_n$ be a sequence of positive numbers. Then, for any $x>0$,
\begin{align}
  {\mathbb P}\Big(  S_n \geq  a_n x\Big)
   &\leq {\mathbb P}\Big(  M_1 \geq  a_n x/2\Big)+ {\mathbb P}\Big(  S_{2,n} \geq  a_n x/2\Big) \nonumber \\
    &\leq {\mathbb P}\bigg(  G_{X_1}(X_1) \geq  \frac{a_n x }{  2 K_{n-1} (\rho)} \bigg) + {\mathbb P}\Big(  S_{2,n} \geq  a_n x/2\Big)   =: I_1(a_n, x)+ I_2(a_n, x) \, , \label{ineq001}
\end{align}
and note that the same bound is valid for ${\mathbb P}(  -S_n \geq  a_n x)$ by replacing the term
$I_2(a_n, x)$ by $\tilde I_2(a_n, x):= {\mathbb P}(  -S_{2,n} \geq  a_n x/2)$.

The  term  $I_1(a_n, x)$ will be most of the time  negligible, and represents  the direct influence of the initial distribution of the chain.
For instance, when the chain starts from a point $X_1=x_1$,  then $G_{X_1}(X_1)=0$ and $I_1(a_n, x)=0.$
The main difficulty is to give an upper bound for  $I_2(a_n, x)$, which is the purpose of the present paper.

\subsection{A first exponential bound}
Under a sub-Gaussian  type condition, we obtain the following proposition.
\begin{prop} \label{BerProp}
Assume that there exists a positive constant $\epsilon$   such that,
for any integer $k\geq 2$,
\begin{equation}\label{Bercond02}
\mathbb{E}\left[ \big(H_{k}(X_{k-1}, \varepsilon_k) \big)^l  \right] \leq \frac{1}{2} \frac{  l! \,\epsilon ^{l-2}  }{(l-1)^{l/2}}  \mathbb{E}\left[  \big( H_{k}(X_{k-1}, \varepsilon_k) \big)^2 \right]
\textrm{ for all}\ l\geq  2.
\end{equation}
 Then,  for any   $ x>0$,
\begin{align}
  \mathbb{P}\Big( \pm S_n \geq x V_n  \Big)
   &\leq I_1(V_n, x)+ \exp\bigg\{ - \frac{(x/2)^2}{  1+ \sqrt{ 1+ x \epsilon  K_{n-2} (\rho)
  /\sigma_n   }  + x   \epsilon K_{n-2} (\rho)/2\sigma_n   } \bigg\} \label{jnsk01}\\
 &\leq  I_1(V_n, x)+ \exp\bigg\{ - \frac{(x/2)^2}{  2 \big( 1+ x   \epsilon K_{n-2} (\rho)  / 2\sigma_n   \big)} \bigg\},\label{jnsk02}
\end{align}
where $$V_n^2 =  \sum_{k=2}^{n} K_{n-k}^2(\rho) \mathbb{E}\big[ \big( H_{k}(X_{k-1}, \varepsilon_k) \big)^2 \big] \ \   and \ \   \sigma_n^2=  \frac1n V_n^2.$$
\end{prop}

\begin{rem} Let us give some comments on Proposition \ref{BerProp}.
\begin{enumerate}
\item Condition \eqref{Bercond02} is in fact a sub-Gaussian condition. On can check that it is satisfied provided
$$   \inf_{k \geq 2} {\mathbb E}\big[ \big( H_{k}(X_{k-1}, \varepsilon_k) \big)^2 \big]    >0  $$
and
$$ \sup_{k \geq 2}  {\mathbb E}\Big[ \exp  \Big\{ c \, \big( H_{k}(X_{k-1}, \varepsilon_k)\big)^2 \Big\} \Big]  < \infty $$ for some positive constant $c$ not depending on $k$.

\item Assume that $${\mathbb E}\Big[\exp  \left\{ c \, \sqrt{G_{X_1}(X_1)} \right\}\Big] < \infty$$ for some positive constant $c$, and that
$$ 0< \liminf_{n\rightarrow \infty}\sigma_n    \leq \limsup_{n\rightarrow \infty}\sigma_n < \infty.$$   Then, it follows from Proposition \ref{BerProp} that
\begin{eqnarray}
\mathbb{P}\left( \pm S_n \geq  n  \right)   = O \left( \exp\left\{ - C \sqrt n  \right\} \right)
\end{eqnarray}
for some positive constant $C$.
\end{enumerate}
\end{rem}

\noindent\emph{Proof.}
 By Taylor's expansion of $e^x$ and the fact that $\mathbb{E} [S_{2,n}]=0$, we have, for all $t\geq 0,$
 \begin{eqnarray}\label{sxsa}
 \mathbb{E}\bigg[\exp\bigg\{t \frac{ S_{2,n}}{\sqrt{n}} \bigg \}\bigg] = 1+ \sum_{k=2}^{\infty} \frac{t^k}{k!}  \mathbb{E}\Big[\Big(\frac{ S_{2,n}}{\sqrt{n}}\Big)^k \Big].
\end{eqnarray}
Using   Rio's inequality (see Theorem 2.1 of \cite{R09}): for any $p\geq 2,$
\begin{eqnarray}\label{R}
\Big(\mathbb{E}[|S_{2,n}|^p] \Big)^{2/p} \leq  (p-1)  \sum_{i=2}^n  \big(\mathbb{E}[|M_i|^p] \big)^{2/p}  ,
\end{eqnarray}
we get, for all $k\geq 2,$
\begin{eqnarray}\label{ssxsfsd}
\mathbb{E}[|S_{2,n}|^k] \leq (k-1)^{k/2} \Big(\sum_{i=2}^n \big(\mathbb{E}[|M_i|^k] \big)^{2/k}\Big)^{k/2}.
\end{eqnarray}
Hence, by H\"{o}lder's  inequality, inequality (\ref{ssxsfsd}) implies that, for all $k\geq 2,$
\begin{eqnarray}\label{sddd}
\mathbb{E}[|S_{2,n}|^k] \leq (k-1)^{k/2} n^{k/2-1} \sum_{i=2}^n  \mathbb{E}[|M_i|^k]  .
\end{eqnarray}
Applying the last inequality to (\ref{sxsa}), we obtain
 \begin{eqnarray}\label{bnlts}
\mathbb{E}\bigg[\exp\bigg\{t \frac{ S_{2,n}}{\sqrt{n}} \bigg \}\bigg] \leq 1+ \sum_{k=2}^{\infty} \Big( \frac{t^k}{k!}  (k-1)^{k/2} n^{ -1 } \sum_{i=2}^n  \mathbb{E}[|M_i|^k] \Big).
\end{eqnarray}
 By points 2 of Proposition \ref{McD} and (\ref{Bercond02}), we deduce that,
for any integer $i\geq 2$,
\begin{eqnarray*}
\mathbb{E}[|M_i|^l] &\leq&  \mathbb{E}[|K_{n-i}(\rho) H_{i}(X_{i-1}, \varepsilon_i)|^l ] \\
&\leq& \frac{1}{2} \frac{  l! \,(K_{n-2}(\rho) \epsilon )^{l-2}  }{(l-1)^{l/2}}   \mathbb{E}[ (K_{n-i}(\rho)H_{i}(X_{i-1}, \varepsilon_i))^2]  \ \ \ \ \  \ \ \textrm{for all}\ l\geq  2.
\end{eqnarray*}
Hence  condition
(\ref{Bercond02}) implies that, for all $0\leq t < (K_{n-2}(\rho) \epsilon)^{-1}$,
 \begin{eqnarray}\label{ineq23}
\mathbb{E}\bigg[\exp\bigg\{t \frac{ S_{2,n}}{\sqrt{n}} \bigg \}\bigg] \leq 1+ \sum_{k=2}^{\infty}  \frac{ \ \sigma_n^2}{2 } \, t^k (K_{n-2}(\rho) \epsilon )^{k-2} = 1+   \frac{t^2 \sigma_n^2  }{2\,(1- t K_{n-2}(\rho) \epsilon )} .
\end{eqnarray}
By the inequality $1+x \leq e^x  ,$ it follows that, for all $0\leq t < (K_{n-2}(\rho) \epsilon)^{-1}$,
 \begin{eqnarray*}
\mathbb{E}\bigg[\exp\bigg\{t \frac{ S_{2,n}}{\sqrt{n}} \bigg \}\bigg] \ \leq \   \exp \Bigg\{ \frac{t^2  \sigma_n^2  }{2\,(1- t  K_{n-2}(\rho)\epsilon)}   \Bigg\}.
\end{eqnarray*}
Applying Markov's inequality,   it is then easy to see that, for all $0\leq t < \sigma_n(K_{n-2}(\rho) \epsilon)^{-1}$ and  $ x \geq 0,$
\begin{eqnarray*}
  \mathbb{P}\left(   S_{2,n} \geq x  V_n /2   \right)
   \leq   \exp\Big\{ - t x /2\Big \} \mathbb{E}\left[   \exp\bigg\{t \frac{ S_{2,n}}{V_n } \bigg \} \right].
\end{eqnarray*}
Hence
\begin{eqnarray*}
  \mathbb{P}\left(   S_{2,n} \geq x  V_n /2   \right)
    &\leq& \inf_{0\leq t <\sigma_n (K_{n-2}(\rho) \epsilon)^{-1}}  \exp\bigg\{ - t x /2 + \frac{  t^2    }{2\,(1- t \, K_{n-2}(\rho)\epsilon/ \sigma_n)}    \bigg \} \\
    &=& \exp\bigg\{ - \frac{(x/2)^2}{  1+ \sqrt{ 1+ x K_{n-2}(\rho)\epsilon
  /\sigma_n   }  + x   K_{n-2}(\rho)\epsilon /2\sigma_n   } \bigg\},
\end{eqnarray*}
which gives   (\ref{jnsk01}).   Using the   inequality
$   \sqrt{ 1+ x K_{n-2}(\rho)\epsilon / \sigma_n    }  \leq  1+  x K_{n-2}(\rho)\epsilon /2 \sigma_n ,$
we get  (\ref{jnsk02})  from (\ref{jnsk01}).  \qed

\subsection{Semi-exponential  bounds}\label{Bersec}
In the case where the variables $H_{k}(X_{k-1}, \varepsilon_k)$ have semi-exponential moments, the following proposition holds.

\begin{prop}\label{propExp}
Let $\alpha \in (0,1)$.
Assume that there exists a constant $C_1$ such that, for any integer $k\geq 2$,
\begin{equation} \label{Bernsteinmoment}
  \mathbb{E}\,\Big[\exp\Big\{ \Big(H_{k}(X_{k-1}, \varepsilon_k)\Big)^{\frac{2\alpha}{1-\alpha}} \Big \}  \Big] \leq C_1  .
\end{equation}
Then, for any $x> 0$,
\begin{eqnarray}\label{mainineq01}
  {\mathbb P}\left(\pm  S_n \geq  n x\right)
   \leq I_1(n, x) +  C(\alpha, x)  \exp\left\{-\left(\frac{ x  }{ 8 K_{n-2}(\rho)}\right)^{2\alpha} n^\alpha\ \right\} \, ,
\end{eqnarray}
where
\[
C(\alpha,  x)=  2+ 35  C_1     \left( \frac{K_{n-2}^{2\alpha}(\rho) }
{ x^{2\alpha} 4^{2-3\alpha}}  + \frac{4 K_{n-2}^2(\rho)}{  x^2 }
\left( \frac{3(1-\alpha)}{2\alpha}\right)^{\frac{ 1-\alpha}{\alpha}} \right)
\]
depends on $n$ only through the term $K_{n-2}(\rho)$.
\end{prop}

\begin{rem}\label{rem01}
Let us comment on inequality  \eref{mainineq01}.
Assume  moreover that
\begin{equation}\label{Bernsteinmoment1}
 \mathbb{E}\,\left[\exp\left\{ c\Big( G_{X_1}(X_1)\Big)^{\alpha  }  \right \} \right]  \leq   C_2 ,
\end{equation}
for two positive constants $c, C_2$.
Then, it follows from \eqref{mainineq01} that
\begin{eqnarray}\label{LU00}
\mathbb{P}\left( \pm S_n \geq  n  \right)   = O \left( \exp\left\{ - C  n^{\alpha}     \right\} \right),
\end{eqnarray}
for some positive constant $C$.
This rate is in accordance with the best possible rate for large deviation of partial sums of martingales differences,
as proved in  Theorem 2.1 of  \cite{Fx1}.
For partial sums of independent random variables, the rate  \eqref{LU00}  holds under weaker conditions on  exponential moments, see Lanzinger and Stadtm\"{u}ller \cite{LS00}.
\end{rem}

\noindent\emph{Proof.}   From point 2 of Proposition \ref{McD} and condition (\ref{Bernsteinmoment}), it is  easy to see that,
   for any  $k \in [2, n],$
\begin{eqnarray}\label{fin36}
\mathbb{E}\,\left[ \exp\left\{  |K_{n-2}^{-1}(\rho)M_k|^{\frac{2\alpha}{1-\alpha}} \right\} \right]  \leq   C_1 .
\end{eqnarray}
Applying Theorem 2.1 of  Fan et al. \cite{Fx1} to the martingale sequence
 $(K_{n-2}^{-1} (\rho)M_k, \mathcal{F}_k)_{k=2,..,n}$, we get, for any $x> 0$,
\begin{eqnarray}
 I_2(n, x)
   \leq C(\alpha, x)  \exp\left\{-\left(\frac{ x  }{ 8 K_{n-2}(\rho)}\right)^{2\alpha} n^\alpha\ \right\} . \label{ineq003}
\end{eqnarray}
Combining the inequalities (\ref{ineq001})  and (\ref{ineq003}),  we obtain  the desired inequality.
   \hfill\qed

 \medskip

   For the next proposition, let us introduce the random variables
 \begin{equation}\label{M2cond}
 L_k(X_{k-1}),  \quad \text{where} \quad L_k(x)= \int (H_k(x, y))^2 P_{\varepsilon_k}(dy) \, ,
 \end{equation}
 and note that
 $$
  L_k(X_{k-1})   = \mathbb{E}\,\left[ \left( H_{k}(X_{k-1}, \varepsilon_k) \right)^2     \big | X_{k-1} \right].
 $$
 According to Proposition \ref{McD}, for any $k \geq 2$, ${\mathbb E}[M_k^2|{\mathcal F}_{k-1}] \leq K_{n-k}(\rho) L_k(X_{k-1})$.

\begin{prop} \label{pr43} Assume  (\ref{c2}), and
let $\alpha \in (0,1)$.
Assume that there exist two constants  $C_1$ and  $C_2$ such that, for any integer $k\geq 2, n \geq 2$,
\begin{equation}\label{cond3.5}
\mathbb{E}\left[\exp\bigg\{ \Big(\frac{1}{n} \sum_{k=2}^n   L_k(X_{k-1})  \Big)^{\frac{\alpha}{1-\alpha} } \bigg \} \right ] \leq C_1
\end{equation}
and
\begin{equation}\label{cond3.6}
  \mathbb{E}\,\left[\exp\left\{ \Big(  H_{k}(X_{k-1}, \varepsilon_k) \Big)^{\frac{\alpha}{1-\alpha}} \right \}  \right] \leq C_2.
\end{equation}
Then, for all  $x> 0$,
\begin{multline}\label{mainineq02}
 \mathbb{P}\left( \pm S_{n}\geq n x  \right)  \leq  I_1(n, x) +
 \exp\left\{-\frac{( x K_{n-2}^{-1}(\rho)  /2)^{1+ \alpha} }{2\big( 1 + x K_{n-2}^{-1}(\rho) /6 \big)} n^\alpha  \right\} \\ +
 (C_1 + nC_2)  \exp\left\{ -  \left(   x K_{n-2}^{-1} (\rho)/2 \right)^\alpha n^\alpha  \right \}   .
\end{multline}
\end{prop}
\begin{rem}
According to Remark \ref{rem01},  under the conditions  \eref{cond3.5},  \eref{cond3.6} and \eref{Bernsteinmoment1},  we have
\begin{eqnarray}
\mathbb{P}\left( \pm S_n \geq  n  \right)   = O \left( \exp\left\{ - C  n^{\alpha}     \right\} \right),
\end{eqnarray}
for some positive constant $C$. This rate is in accordance with the best possible rate for large deviation of partial sums of martingales differences,
as proved in  Corollary 2.3  of  \cite{Fx1}.
Note that if $\alpha \in [1/2,1)$, the condition \eqref{cond3.5} is true provided that
\begin{equation}\label{bisbis}
 \sup_{k \geq 2} \mathbb{E}\left[\exp\bigg\{  ( L_k(X_{k-1}) )^{\frac{\alpha}{1-\alpha} } \bigg \} \right ] \leq C_1 \, .
\end{equation}
For $\alpha \in [1/2, 1)$, the two conditions \eqref{bisbis} and \eqref{cond3.6} are clearly less restrictive than \eqref{Bernsteinmoment}, so   Proposition \ref{pr43}
is more precise than Proposition \ref{propExp} in the regime of large deviation. However, it does not allow to control moderate deviations
$\mathbb{P}\left( \pm S_n \geq  n^\delta  \right)$ for $\delta$ close to $1/2$, which is possible via Proposition \ref{propExp} (see
for instance \eqref{etlabete}).
\end{rem}

\medskip

\noindent\emph{Proof.} From point 2 of Proposition \ref{McD} and condition (\ref{cond3.6}), it is  easy to see that, for any  $k \in [2, n],$
\begin{eqnarray*}
\mathbb{E}\,\left[ \exp\left\{  |K_{n-2}^{-1}(\rho)M_k|^{\frac{\alpha}{1-\alpha}} \right\} \right]
 \leq  C_2.
\end{eqnarray*}
For any  $k \in [2, n],$
\begin{eqnarray}
\mathbb{E}\,\left[ |K_{n-2}^{-1}(\rho) M_k|^2 \Big| \mathcal{F}_{k-1} \right]    &\leq&   \mathbb{E}\,\left[ \left(K_{n-2}^{-1}(\rho) K_{n-k}(\rho)H_{k}(X_{k-1}, \varepsilon_k) \right)^2
 \Big| \mathcal{F}_{k-1} \right]     \nonumber\\
 &\leq&   \mathbb{E}\,\Big[ \Big( H_{k}(X_{k-1}, \varepsilon_k) \Big)^2     \Big| \mathcal{F}_{k-1} \Big]  .
\end{eqnarray}
Thus
$$\sum_{k=2}^n \mathbb{E}\,\Big[ |K_{n-2}^{-1}(\rho) M_k|^2 \Big| \mathcal{F}_{k-1} \Big] \leq    \sum_{k=2}^n \mathbb{E}\,\Big[ \Big( H_{k}(X_{k-1}, \varepsilon_k) \Big)^2     \Big| \mathcal{F}_{k-1} \Big] =  \sum_{k=2}^n L_k(X_{k-1}) .$$
Using Theorem 2.2 of Fan et al. \cite{Fx1}, we have
\begin{multline*}
  \mathbb{P}\bigg( K_{n-2}^{-1}(\rho) S_{2,n}\geq n K_{n-2}^{-1}(\rho) x/2\ \ \mbox{and}\ \ \sum_{k=2}^n L_k(X_{k-1}) \leq n
 v^2\    \bigg) \\
   \leq \exp\left\{-\frac{(  K_{n-2}^{-1}(\rho) x/2)^2 }{2\big(   n^{\alpha-1 }v^2+\frac{1}{3}(  K_{n-2}^{-1}(\rho) x/2)^{2-\alpha} \big)} n^\alpha  \right\}  + nC_2  \exp\Bigg\{ -  \Big( \frac{   x }{2 K_{n-2}(\rho)} \Big)^\alpha n^\alpha  \Bigg \}.
\end{multline*}
From the last inequality, we deduce that
\begin{eqnarray*}
I_2(n, x)\!\!  &\leq& \!\! \exp\left\{-\frac{(  K_{n-2}^{-1}(\rho) x/2)^2 }{2\big(   n^{\alpha-1 }v^2+\frac{1}{3}(  K_{n-2}^{-1}(\rho) x/2)^{2-\alpha} \big)} n^\alpha  \right\}  + nC_2  \exp\Bigg\{ -  \Big( \frac{   x }{2 K_{n-2}(\rho)} \Big)^\alpha n^\alpha  \Bigg \}  \\
 &&\, + \, \mathbb{P}\bigg(   \sum_{k=2}^n L_k(X_{k-1})  >
 nv^2\    \bigg).
\end{eqnarray*}
Using the exponential Markov  inequality and the condition \eqref{cond3.5},
we get, for all $ v>0$,
\begin{eqnarray}
\mathbb{P}\bigg(   \sum_{k=2}^n L_k(X_{k-1})  >
 n v^2\    \bigg)\leq C_1 \exp\Big\{ -  v^{\frac{2\alpha}{1-\alpha} } \Big \}. \nonumber
\end{eqnarray}
Taking $v^2 = (x n  K_{n-2}^{-1} (\rho)/2)^{(1-\alpha)/2}$, we have, for all $x>0$,
$$
I_2(n, x) \leq
\exp\left\{-\frac{( x K_{n-2}^{-1}(\rho)  /2)^{1+ \alpha} }{2\big( 1 + x K_{n-2}^{-1}(\rho) /6 \big)} n^\alpha  \right\}  +
 (C_1 + nC_2)  \exp\left\{ -  \left(   x K_{n-2}^{-1} (\rho)/2 \right)^\alpha n^\alpha  \right \}  \, .
$$
Combining the last inequality and   (\ref{ineq001}), we obtain the desired  inequality. \hfill\qed

\subsection{Fuk-Nagaev type  bound}
We now consider the case where the random variables
$H_{k}(X_{k-1}, \varepsilon_k), k\geq 2,$ have only a weak moment of order $p > 2$.
For any real-valued random variable $Z$ and any $p\geq 1$, define the weak
moment of order $p$  by
\begin{equation}\label{weakp}
\|Z\|_{w,p}^p=\sup_{x>0} x^p{\mathbb P}(|Z|>x)\, .
\end{equation}
The following proposition is a  Fuk-Nagaev type inequality (cf.\  Fuk \cite{F73} and Nagaev \cite{N79}; see also Fan \textit{et al.}\ \cite{FGL17} and Rio \cite{R17} for martingales).
\begin{prop}\label{pro3}
Let  $p\geq 2$ and $\delta>0$, and consider the variables $L_k(X_{k-1})$  defined in \eqref{M2cond}.
Assume that there exist two constants  $C_1$ and  $C_2$ such that, for any integer $k\geq 2$,
\begin{equation}\label{cond3.7}
\left \|  \frac{1}{n} \sum_{k=2}^n L_k(X_{k-1})  \right \|_{w,p+ \delta}^{p+ \delta} \leq C_1
\end{equation}
and
\begin{equation}\label{cond3.8}
\left \| H_{k}(X_{k-1}, \varepsilon_k) \right \|_{w,p+ \delta}^{p+ \delta} \leq C_2.
\end{equation}
Then, for all  $x> 0$,
\begin{equation}
 \mathbb{P}\left( \pm S_{n}\geq n x  \right)    \leq I_1(n, x) +  \exp\left\{-\frac{(   K_{n-2}^{-1}(\rho)  /2)^2 }{2\left(n^{-1/(p+\delta)}  x^{-1}+\frac{1}{6}   K_{n-2}^{-1}(\rho)    \right)} (nx)^{\delta/(p+\delta) }  \right\}
    + \  \frac{C_1+ C_2}{n^{p-1} x^{p }}  \, .   \label{sfdw}
\end{equation}
\end{prop}

\begin{rem}\label{commentFNI}
 Let us comment on Proposition \ref{pro3}.
\begin{enumerate}
\item If  there exists a constant $C_3$ such that
\begin{equation}\label{seinsntq}
\left|\left|G_{X_1}(X_1) \right|\right|_{w,p-1}^{p-1} \leq   C_3 ,
\end{equation}
then, for any $x> 0$,
\begin{eqnarray}
 I_1(n, x) ={\mathbb P}\bigg(  G_{X_1}(X_1) \geq  \frac{n x }{  2 K_{n-1} (\rho)} \bigg)\leq (2 K_{n-1} (\rho))^{p-1}  \frac{C_3}{(nx)^{p-1}}  .
\end{eqnarray}
 Thus under  conditions \eref{cond3.7},  \eref{cond3.8} and  \eref{seinsntq},   we have
\begin{equation*}
\mathbb{P}\left( \pm S_n \geq  n  \right)   = O \left( \frac{1}{n^{p-1}}  \right)   .
\end{equation*}
 \item Assume  moreover that $F_n$ satisfies \eref{c2}. Then, according to Remark \ref{remc2},
 Proposition \ref{pro3} remains valid when   $H_{k}(X_{k-1}, \varepsilon_k)$ is replaced by $C(X_{k-1})  G_k(\varepsilon_k) .$
 Since  $C(X_{k-1})$ and $G_k(\varepsilon_k)$ are independent, we easily see that
 $$ L_k(X_{k-1}) \leq (C(X_{k-1}) )^2\mathbb{E} \left [ (G_k(\varepsilon_k))^2\right ] $$
 and
 \begin{equation}
  \left \| C(X_{k-1})  G_k(\varepsilon_k) \right \|^{p+\delta}_{w, p + \delta}    \leq {\mathbb E} \left [ \big(C(X_{k-1}) \big)^{p+\delta} \right]
  \left \| G_k(\varepsilon_k)\right \|^{p+\delta}_{w, p + \delta} \,  .
\end{equation}
Thus if
 \begin{equation*}\label{trivial1}
 {\mathbb E} \left [ \big( C (X_{k-1}) \big) ^{p+\delta}    \right]  \leq C_3  \ \ \ \textrm{and} \ \ \  \left \| G_k(\varepsilon_k)\right \|^{p+\delta}_{w, p + \delta} \leq C_4 \, ,
\end{equation*}
then  condition (\ref{cond3.8}) is satisfied with $C_2=C_3C_4.$ Of course, the same computations may be done by interchanging $C(X_{k-1})$
and $G_k(\varepsilon_k)$. Hence, if
\begin{equation*}\label{trivial2}
 {\mathbb E} \left [ \big( G_k(\varepsilon_k) \big)^{p+\delta}    \right]  \leq C_3  \ \ \ \textrm{and} \ \ \  \left \| C(X_{k-1})\right \|^{p+\delta}_{w, p + \delta} \leq C_4 \, ,
\end{equation*}
then  condition (\ref{cond3.8}) still holds with $C_2=C_3C_4.$
\end{enumerate}
\end{rem}

\medskip

\noindent\emph{Proof.}
To prove Proposition \ref{pro3}, we need the following inequality whose proof can be found in Fan  \textit{et al.}\  \cite{FGL12} (cf. Corollary 2.3 and Remark 2.1 therein).
\begin{lem}
\label{leedm1} Assume that $(\xi_i, {\mathcal G}_i)_{i\geq 1}$ are square integrable martingale differences, and let
$Z_n= \xi_1 + \cdots + \xi_n$ and $ \langle Z \rangle_n = \sum_{k=1}^n {\mathbb E}[\xi_k^2 |{\mathcal G}_{k-1}]$.   Then,
for all $ x, y, v > 0$,
\begin{equation*}
{\mathbb P}\Big( Z_n \geq x  \  \textrm{and} \ \langle Z\rangle_n \leq v^2   \Big)
  \leq    \exp\left\{- \frac{x^2}{2(v^2+ \frac{1}{3}xy)} \right\}+ \mathbb{P}\left( \max_{1\leq i \leq n} \xi_{i}>y  \right).
\end{equation*}
\end{lem}
\noindent By Lemma \ref{leedm1}  and Markov's inequality,  it follows that, for all $ x, y, v > 0$,
\begin{multline*}
  \mathbb{P}\left( K_{n-2}^{-1}(\rho) S_{2,n}\geq n K_{n-2}^{-1}(\rho) x/2\ \ \mbox{and}\ \ \sum_{k=2}^n L_k(X_{k-1}) \leq n v^2\    \right) \\
   \leq \exp\left\{-\frac{( n K_{n-2}^{-1}(\rho) x/2)^2 }{2\big(  n v^2+\frac{1}{6}  n K_{n-2}^{-1}(\rho) x y   \big)}   \right\}  +\mathbb{P}\left( \max_{2\leq i \leq n} K_{n-2}^{-1}(\rho) M_i \geq y  \right).
\end{multline*}
It is easy to see that,  for all $   y  > 0$,
\begin{align*}
\mathbb{P}\left( \max_{2\leq i \leq n} K_{n-2}^{-1}(\rho) M_i \geq y  \right)
 & \leq  \sum_{i=2}^n \mathbb{P}\left(  K_{n-2}^{-1}(\rho) M_i \geq y  \right)
  \leq  n  C_2 y^{-(p+\delta)} ,
\end{align*}
and that, for all $v>0$,
\begin{equation*}
  \mathbb{P}\left( \sum_{k=2}^n L_k(X_{k-1})    > n v^2   \right)
   \leq  C_1  v^{-2(p+\delta)}.
\end{equation*}
Thus,  for all $ x, y, v > 0$,
\begin{eqnarray*}
I_2(n, x)  \leq
   \exp\left\{-\frac{( n K_{n-2}^{-1}(\rho) x/2)^2 }{2\big(   n v^2+\frac{1}{6}  n K_{n-2}^{-1}(\rho) x y   \big)}   \right\}   + \frac{nC_2}{ y^{p+\delta}}  +  \frac{C_1}{ v^{2(p+\delta)}  } .
\end{eqnarray*}
Taking $y =   (nx)^{p/(p+\delta) }$ and $v^{2(p+\delta)} = n^{p-1} x^p $  in the last inequality, we obtain,
 for all $ x > 0$,
\begin{eqnarray*}
I_2(n, x)  \leq
   \exp\left\{-\frac{(   K_{n-2}^{-1}(\rho)  /2)^2 }{2\big( n^{-1/(p+\delta)}  x^{-1} +\frac{1}{6}   K_{n-2}^{-1}(\rho)    \big)} (nx)^{\delta/(p+\delta) }  \right\}   +   \frac{C_1+ C_2}{n^{p-1} x^{p }}   .
\end{eqnarray*}
Combining the last inequality and   (\ref{ineq001}), we obtain the desired  inequality. \hfill\qed

\subsection{von Bahr-Esseen's inequality, weak form}
We now consider the case where the variables $G_{X_1}(X_1)$ and
$H_k(X_{k-1}, \varepsilon_k)$ have only a weak moment of order $p \in (1,2)$.

\begin{prop}\label{weakVBEI}
Let $p \in (1, 2)$. Assume that there exists positive constants $A_k(p)$ such that, for any $k \in [2,n]$,
\begin{equation}\label{weakVB}
  \left \|H_k(X_{k-1}, \varepsilon_k)\right \|_{w,p}^p \leq A_k(p)
  \, .
\end{equation}
Then, for any $x>0,$
\begin{equation}\label{weakVBE}
 {\mathbb P}(|S_n| \geq x) \leq 2 I_1(1,x)+ \frac{2^pC_p B_p(n,\rho)}{x^p}\, ,
\end{equation}
where
 $$
 C_p= \frac{4p}{p-1} +\frac{8}{2-p}\,
$$
and
 $$
B_p(n,\rho)=   \sum_{k=2}^n \left(K_{n-k}(\rho)\right)^p  A_k(p) \,.
$$
\end{prop}

\begin{rem}
Let us comment on Proposition \ref{weakVBEI}.
\begin{enumerate}
\item Contrary to the previous inequalities of Section \ref{deviationiq}, this inequality is truly non-stationary, in the sense that
it is expressed in terms of the weak moments $\|H_k(X_{k-1}, \varepsilon_k) \|_{w,p}^p$, without assuming a uniform bound (in $k$) on these moments.
This will also be the case of the moment inequalities of Section \ref{SECMOM}.
\item Assume  moreover that $F_n$ satisfies \eref{c2}. Then, according to Remark \ref{remc2},
 Proposition \ref{pro3} remains true if  condition   \eref{weakVB} is replaced by
\begin{equation}\label{weakVfdB}
  \left \| C(X_{k-1})  G_k(\varepsilon_k) \right \|_{w,p}^p \leq A_k(p) \, .
\end{equation}
In particular, if either
$$
{\mathbb E} \left [  (C(X_{k-1}))^p       \right]  \leq A_{1,k}(p)  \ \ \ \textrm{and} \ \ \  \left \| G_k(\varepsilon_k)\right \|^{p}_{w, p} \leq A_{2,k}(p)
$$
or
$$
 {\mathbb E} \left [  (G_k(\varepsilon_k))^p   \right]  \leq A_{1,k}(p)  \ \ \ \textrm{and} \ \ \  \left \| C(X_{k-1})\right \|^{p}_{w, p} \leq A_{2,k}(p) \, ,
$$  hold,
then condition \eqref{weakVfdB} is satisfied  with $A_k(p)=A_{1,k}(p)A_{2,k}(p)$.

\item Assume that  $\|   G_{X_1}(X_1) \|_{w,p-1}< \infty$ and that
$B_p(n,\rho) =O(n)$, then
\begin{equation*}
\mathbb{P}\left( |S_n| \geq  n  \right)   = O \left( \frac{1}{n^{p-1}}  \right)  .
\end{equation*}

\end{enumerate}
\end{rem}

\noindent\emph{Proof.} By Proposition 3.3 of Cuny \emph{et al.}  \cite{CDM17}, we have, for  any $x>0,$
\begin{equation}\label{inew01}
 {\mathbb P}(|S_{2,n}|\geq x/2) \leq  \frac{2^p C_p}{x^p}\, \sum_{k=2}^n   \left\| M_k\right\|_{w,p}^p  .
\end{equation}
From point 2 of Proposition \ref{McD} and condition (\ref{weakVB}),  it follows that,
for any  $k \in [2, n],$
\begin{equation}\label{inew03}
 \left \| M_k\right\|_{w,p}^p  \leq   \left \|  K_{n-k}(\rho)     H_{k}(X_{k-1}, \varepsilon_k) \right \|_{w,p}^p  \leq   \left(K_{n-k}(\rho) \right)^p   A_k(p).
\end{equation}
Combining  \eref{inew01}  and  \eref{inew03}, we obtain the desired inequality.  \hfill\qed

\section{Moment inequalities}\label{SECMOM}

In this section, we shall control the ${\mathbb L}_p$-norm of the functional $S_n$, for $p>1$. The upper bounds will be expressed in terms of
the moments of order $p$ of the dominating random variables $H_{k} (X_{k-1}, \varepsilon_k)$. Let us emphasize that all the inequalities of these
section are completely non-stationary, in the sense that we shall  not assume a uniform bound (in $k$) on $\|H_k(X_{k-1}, \varepsilon_k) \|_{p}^p$.

\subsection{Marcinkiewicz-Zygmund type inequality}
We assume in this subsection  that the dominating random variables $G_{X_1}(X_1) $  and $H_{k} (X_{k-1}, \varepsilon_k)$ have a moment of order $p \geq 2$.
\begin{prop}\label{MZIP}
Let $p \geq 2$. Assume that there exist    positive constants $A_k(p)$   such that
\begin{equation}\label{shnbs}
 \mathbb{E}  \Big[\big(G_{X_1}(X_1)\big)^{p} \Big]  \leq A_1(p), \quad \text{and for $k \in [2,n]$} \quad \mathbb{E} \Big[ \big(H_{k} (X_{k-1}, \varepsilon_k)\big)^p  \Big]   \leq A_k(p)  .
\end{equation}
Then
\begin{equation}\label{MZI}
  \| S_n \|_p \leq   \sqrt{ A_p(n,\rho)}\, ,
\end{equation}
 where
 $$
A_p(n,\rho)= \big(K_{n-1}(\rho) \big)^2 \big(A_1(p) \big)^{2/p} +
(p-1) \sum_{k=2}^n \big(K_{n-k}(\rho) \big)^2  \big(A_k(p) \big)^{2/p}\, .
$$
\end{prop}
\begin{rem}
Assume  moreover that $F_n$ satisfies \eref{c2}. Then, according to Remark \ref{remc2},
 inequality \eref{MZI} remains true if the second condition of  \eref{shnbs} is replaced by
\begin{equation} \label{weaker}
  {\mathbb E} \left [(C(X_{k-1}))^p \right ] {\mathbb E} \left [ (G_k(\varepsilon_k))^p \right ]
  \leq A_k(p) \, .
\end{equation}
\end{rem}

\noindent\emph{Proof.} Applying Theorem 2.1 of Rio \cite{R09}, we get
 $$
\| S_n \|_p^2  \leq \|M_1\|_p^2 +(p-1) \sum_{k=2}^n  \| M_k \|_p^2 \, .
$$
By Proposition \ref{McD}  and condition (\ref{shnbs}), it follows that
\begin{eqnarray*}
\| S_n \|_p^2   &\leq&   \big(K_{n-1}(\rho)\big)^2  \left( {\mathbb E} \left[  \big(  G_{X_1}(X_1)\big)^p\right]  \right)^{\frac 2 p} + (p-1)
    \sum_{k=2}^n \big(K_{n-k}(\rho)\big)^2 \big({\mathbb E} \left[  \big( H_{k}(X_{k-1}, \varepsilon_k)   \big)^p\right]\big)^{\frac 2 p} \\
    &\leq&   \big(K_{n-1}(\rho)\big)^2  \big( {\mathbb E}  \left[  \big(  G_{X_1}(X_1)\big)^p\right]  \big)^{\frac 2 p} + (p-1)
    \sum_{k=2}^n \big(K_{n-k}(\rho)\big)^2 \big( {\mathbb E} \left[ \big(H_{k}(X_{k-1}, \varepsilon_k)\big)^p \right]\big)^{\frac 2 p} \\
    &\leq&  A_p(n,\rho),
\end{eqnarray*}
which gives the desired inequality. \hfill\qed

\subsection{Rosenthal's inequality}
Under the same assumptions as in the previous subsection, one can prove the following Rosenthal-type inequality.
\begin{prop} \label{Ros}
Let $p\geq 2$, and consider the variables $L_k(X_{k-1})$  defined in \eqref{M2cond}. If  \eqref{shnbs} holds,
then there exists a constant $C_p$ depending only on $p$ such that
\begin{equation}
\|S_n\|^p_p   \leq C_p \bigg(  \mathbb{E}\Big[ \Big( \big(K_{n-1}(\rho) \big)^2 A_1(2)  +     \sum_{k=2}^n \big( K_{n-k} (\rho)\big)^2 L_k(X_{k-1})  \Big)^{p/2} \Big]
+   \sum_{k=1}^n \big(K_{n-k}(\rho) \big)^p  A_k(p) \bigg).  \label{vBEs}
\end{equation}
\end{prop}
\begin{rem}
Assume that $F_n$ satisfies \eref{c2}.  Then, according to Remark \ref{remc2},
 it follows from the proof of Proposition \ref{Ros} that
 inequality \eref{vBEs} remains true if the second condition of  \eref{shnbs} is replaced by \eqref{weaker}.
\end{rem}
\noindent\emph{Proof.} From point 2 of Proposition \ref{McD},  it is easy to see that
\begin{equation*}
   \sum_{k=1}^n\mathbb{E}[ M_k^2 |\mathcal{F}_{k-1}]  \leq \big( K_{n-1} (\rho)\big)^2 \mathbb{E}  [ \big(G_{X_1}(X_1) \big)^{2} ]+  \sum_{k=2}^n \big( K_{n-k}(\rho) \big)^2
   L_k(X_{k-1}) \, ,
\end{equation*}
and that
\begin{align*}
   \sum_{k=1}^n\mathbb{E}[| M_k|^p]  &\leq \big( K_{n-1} (\rho) \big)^p \mathbb{E} [ \big( G_{X_1}(X_1) \big)^{p} ]+    \sum_{k=2}^n \big( K_{n-k} (\rho)\big)^p \mathbb{E} \big[  \big( H_{k} (X_{k-1}, \varepsilon_k) \big)^p   \big] \\
   & \leq \big( K_{n-1}(\rho)\big)^p A_1(p)  +     \sum_{k=2}^n  \big(K_{n-k}(\rho) \big)^p A_k(p)  .
\end{align*}
The desired inequality is then a direct consequence of
Rosenthal's inequality for martingales (see for instance Theorem 2.12 of Hall and Heyde \cite{HH80}).   \hfill\qed

\subsection{von Bahr-Esseen's inequality}
In this subsection, we assume that the dominating random variables $G_{X_1}(X_1)$
and $H_{k}  (X_{k-1}, \varepsilon_k)$ have  a moment of order $p \in (1,2]$.
\begin{prop}\label{VBEI}
Let  $p \in (1, 2]$,  and assume that  \eqref{shnbs} holds.
Then
\begin{equation}\label{vBE}
  \| S_n \|_p^p  \leq     A_p(n,\rho) ,
\end{equation}
 where
\begin{equation}\label{Ap}
A_p(n,\rho)= \big( K_{n-1}(\rho) \big )^p A_1(p)  +
2^{2-p}   \sum_{k=2}^n \big( K_{n-k}(\rho) \big)^p  A_k(p) \, .
\end{equation}
\end{prop}

\begin{rem}
Assume that $F_n$ satisfies \eref{c2}. Then,   according to Remark \ref{remc2}, it follows from the proof of Proposition \ref{VBEI} that
  inequality \eref{vBE} remains true if the second condition of  \eref{shnbs} is replaced by \eqref{weaker}\end{rem}
\noindent\emph{Proof.} Using an improvement of the von Bahr-Esseen inequality  (see  inequality (1.11) in Pinelis \cite{P10}), we have
$$
\| S_n \|_p^p \leq \|M_1\|_p^p +  \tilde C_p \sum_{k=2}^n  \| M_k \|_p^p \, ,
$$
where the constant $\tilde C_p$ is described in Proposition 1.8 of Pinelis \cite{P10}, and is such that $\tilde C_p \leq 2^{2-p}$
for any $p \in [1,2]$.
By Proposition \ref{McD}, it follows that
\begin{align*}
\| S_n \|_p^p   &\leq    \big(  K_{n-1}(\rho) \big)^p   {\mathbb E} \left[  \big(  G_{X_1}(X_1) \big)^p \right]    + \tilde C_p
    \sum_{k=2}^n \big( K_{n-k}(\rho) \big)^p {\mathbb E} \left[ \big(  H_{k}(X_{k-1}, \varepsilon_k) \big)^p  \right] \\
    &\leq   \big( K_{n-1}(\rho) \big)^p A_1(p)   +
\tilde C_p    \sum_{k=2}^n \big( K_{n-k} (\rho)\big)^p A_k(p) \, ,
\end{align*}
which is the desired inequality. \hfill\qed

\medskip

\paragraph{Aknowledgements} P. Doukhan wants to thank Tianjin  University for its hospitality.
His work has been developed within the MME-DII center of excellence (ANR-11-LABEX-0023-01) \& PAI-CONICYT MEC N$^o$ 80170072.
X. Fan has been partially  supported by the  National Natural Science Foundation of China (Grant no.\ 11601375).

\end{document}